\documentclass[12pt]{article}

\usepackage{amssymb,amsmath,amsthm, mathtools}
\usepackage{amsfonts}
\usepackage{enumerate}

\usepackage{authblk}

\usepackage{palatino}
\usepackage{fancyhdr}	
\usepackage{mdframed}

\usepackage{geometry}
\newtheoremstyle{mytheoremstyle} 
    {10pt}                    
    {10pt}                    
    {\normalfont}                   
    {}                           
    {\bfseries}                   
    {.}                          
    {0.3cm}                       
    {}  
\theoremstyle{mytheoremstyle}

\theoremstyle{plain}

\fancypagestyle{plain}{%
  \fancyhf{}%
}

\newtheorem{theorem}{Theorem}[section]

\newtheorem{lemma}[theorem]{Lemma}
\newtheorem{proposition}[theorem]{Proposition}

\newtheorem{remark}[theorem]{Remark}
\newtheorem{question}[theorem]{Question}

\newtheorem{problem}[theorem]{Problem}
\newtheorem{example}[theorem]{Example}

\begin{document}
\title{On unimodular tournaments}
\author{Wiam Belkouche,
	Abderrahim Boussa\"{\i}ri\thanks{Corresponding author: Abderrahim Boussaïri. Email: aboussairi@hotmail.com} , Abdelhak Chaïchaâ and Soufiane Lakhlifi
}
\affil{Laboratoire Topologie, Alg\`ebre, G\'eom\'etrie et Math\'ematiques Discr\`etes, Facult\'e des Sciences A\"in Chock, Hassan II University of Casablanca, Maroc.}

\maketitle

\begin{abstract}
A tournament is unimodular if the determinant of its skew-adjacency 
matrix is $1$. In this paper, we give some properties and constructions 
of unimodular tournaments. A unimodular tournament $T$ with skew-adjacency
matrix $S$ is invertible if $S^{-1}$ is the skew-adjacency 
matrix of a tournament. A spectral characterization of invertible 
tournaments is given. Lastly, we show that every $n$-tournament can be 
embedded in a unimodular tournament by adding at most $n - \lfloor\log_2(
n)\rfloor$ vertices.
\end{abstract}



\section{Introduction}
  Let $T$ be a tournament with vertex set $\{v_1,\ldots,v_n\}$. The \emph
{skew-adjacency matrix} of $T$ is the $n\times n$ zero-diagonal matrix	$S
=[s_{ij}]_{1\leq i,j\leq n}$ in which $s_{ij}=1$ and $s_{ji}=-1$ if $v_i$
 dominates $v_j$. Equivalently, $S=A-A^{t}$ where $A$ is the adjacency 
matrix of $T$. We define the determinant $\det(T)$ of $T$ as the 
determinant of $S$. As $S$ is skew symmetric, $\det(T)$ vanishes when $n$ 
is odd. When $n$ is even, the determinant is the square of the Pfaffian 
of $S$. Moreover, McCarthy and Benjamin \cite[Proposition~1]{mccarthy96}
proved that the determinant of an $n$-tournament has the same parity as $n-1$.

\begin{proposition}\label{mccarthy}
	The determinant of a tournament with an even number of vertices is the
	square of an odd number.
\end{proposition}
		
	Let $T$ be a tournament. The \emph{converse} of $T$, obtained by 
reversing all its arcs, has the same determinant as $T$. The switching 
is another operation that preserves the determinant. The \emph{switch} 
of a tournament on a vertex set $V$ with respect to a subset $X$, is the 
tournament obtained by reversing all the arcs between $X$ and $V\setminus
 X$. It is well-known that if two tournaments are switching equivalent, 
then their skew-adjacency matrices are $\{\pm 1\}$-diagonally similar 
\cite{moorhouse95}. Hence, switching equivalent tournaments have the 
same determinant. 

	In this paper, we consider the class of tournaments whose skew-adjacency
matrices are unimodular, or equivalently, tournaments whose 
determinants are equal to one. We call such tournaments \emph{unimodular}.
By the forgoing, this class is closed under the converse and switching 
operations. Examples of unimodular tournaments are transitive 
tournaments with an even number of vertices and their switches. The 
smallest tournaments that are not unimodular consist of a vertex 
dominating or dominated by a $3$-cycle. These tournaments are called 
\emph{diamonds} \cite{gnanvo1992reconstruction}. A tournament contains 
no diamonds if and only if it is switching equivalent to a transitive 
tournament \cite{babai2000automorphisms}. Tournaments without diamonds 
are known as \emph{local orders} \cite{cameron1978orbits}, \emph{locally 
transitive tournaments} \cite{lachlan1984countable}, or \emph{vortex-free
tournaments} \cite{knuth1992axioms}.
		
	The \emph{join} of a tournament $T_1$ to a tournament $T_2$, denoted 
by $T_1 \rightarrow T_2$, is the tournament obtained from $T_1$ and $T_2$
 by adding an arc from each vertex of $T_1$ to all vertices of $T_2$. 
The join of a tournament $T$ to a tournament with one vertex is denoted 
by $T^{+}$. Our first main result gives a necessary and sufficient 
condition on the unimodularity of the join of two tournaments. It follows
directly from Theorem \ref{main7}, which will be proved in the next section.

\begin{theorem}\label{main6}
  Let $T_1$ and $T_2$ be two tournaments with $p$ and $q$ vertices respectively.
  \begin{enumerate}[i)]
    \item If $p$ and $q$ are even, then $T_1 \rightarrow T_2$ is unimodular
		if and only if $T_1$ and $T_2$ are unimodular.
    \item If $p$ and $q$ are odd, then $T_1 \rightarrow T_2$ is unimodular
		if and only if $T_1^{+}$ and $T_2^{+}$ are unimodular.
  \end{enumerate}
\end{theorem}
	
  Let $T$ be a unimodular tournament and let $S$ be its skew-adjacency 
matrix. The inverse $S^{-1}$ of $S$ is a unimodular skew-symmetric 
integral matrix, but its off-diagonal entries are not necessarily from $
\{-1, 1\}$.  We say that $T$ is invertible if $S^{-1}$ is the skew-adjacency
matrix of a tournament. We call this tournament the inverse of 
$T$ and we denote it by $T^{-1}$. For graphs, the inverse was introduced 
by considering the adjacency matrix and has been studied extensively 
\cite{cvetkovic1978self, godsil1985inverses, kirkland2007unimodular, 
barik2006nonsingular, tifenbach2009directed}.
  
	We give a spectral characterization of invertible tournaments. 
Moreover, we prove that every $n$-tournament can be embedded in an 
invertible, and hence unimodular, $2n$-tournament. The following problem 
arises naturally.
  
\begin{problem}
	For a tournament $T$ on $n$ vertices, what is the smallest number $u^{+
}(T)$ of vertices we must add to $T$ to obtain a unimodular tournament?
\end{problem}
  
  We prove that $u^{+}(T)$ cannot exceed $n - \lfloor\log_2(n)\rfloor$. 
Moreover, if the skew-adjacency matrix $S$ of $T$ is a \emph{skew-conference 
matrix}, that is, $S^2 = (1-n)I_n$, then $u^{+}(T)$ is at least $n/2$.
Hence, $u^{+}(T)$ can be arbitrarily large.

\section{The determinant of the join of tournaments}\label{join}
  Let $T_1$ and $T_2$ be two tournaments, and let $\chi_1(x)$ and $\chi_2
(x)$ be the characteristic polynomials of their adjacency matrices. 
Then, the characteristic polynomial of $T_1\rightarrow T_2$ is $\chi_1(x)
\chi_2(x)$. There is no similar result for the skew-adjacency matrix. 
However, we obtain the following result.
 
\begin{theorem}\label{main7}
  Let $T_1$ and $T_2$ be two tournaments with $p$ and $q$ vertices respectively.
  \begin{enumerate}
    \item If $p$ and $q$ are even, then $\det(T_1 \rightarrow T_2) = \det(T_1)\cdot \det(T_2)$.
    \item If $p$ and $q$ are odd, then $\det(T_1 \rightarrow T_2) = \det(T_1^{+})\cdot \det(T_2^{+})$.
  \end{enumerate}
\end{theorem}
	
	Let $S_1$ and $S_2$ be the skew-adjacency matrices of $T_1$ and $T_2$ 
respectively. The skew-adjacency matrix $S$ of $T_1 \rightarrow T_2$ can 
be written as follows

\begin{equation}
  S=
		\left(
					\begin{array}{c|c}
						S_1 & J \\
						\hline \vspace{0.1cm}
						-J^t & S_2
					\end{array}
		\right)\mbox{.}
\end{equation}
where $J$ is the all ones matrix of order $p\times q$. As mentioned 
above, if $p$ is even, $S_1$ is non-singular. The first assertion 
follows from the more general result.

\begin{lemma}\label{lem2}
	Let $M$ be a skew-symmetric matrix of the form
	\[
		M=\begin{pmatrix}
				A & B \\
				-B^t & D 
			\end{pmatrix}\mbox{.}
	\]
If $A$ is non singular and $rank(B)=1$, then \[ \det(M)= \det(A)\cdot\det(D)\mbox{.}\]
\end{lemma}

\begin{proof}
	Using Schur's complement formula, we get \[\det(M)=\det(A)\cdot\det(D+BA^{-1}B^t)\mbox{.}\]
	As $rank(B) = 1$, there exist two column vectors $\alpha$ and $\beta$ 
such that $B = \alpha\beta^t$. Then, we have 
	 $$ BA^{-1}B^t  =  \alpha\beta^t A^{-1} \beta\alpha^t$$
	Since the matrix $A^{-1}$ is skew-symmetric and $(\beta^t A^{-1} \beta)$
is a scalar, $\beta^t A^{-1} \beta=0$, and hence \[\det(M)=\det(A)\cdot\det(D)\mbox{.}\qedhere\]
\end{proof}
	
\begin{proof}[Proof of Theorem \ref{main7}]
  The first assertion is already proven. For the second assertion, we 
prove that $\det(T_1\rightarrow T_2) = \det(T_1^{+}\rightarrow T_2^{+})$.
For this, consider the tournament $R$ on two vertices. It is easy to 
see that $T_1^{+}\rightarrow T_2^{+}$ is switching equivalent to $(T_1
\rightarrow T_2) \rightarrow R$. Hence, $\det(T_1^{+}\rightarrow T_2^{+}
)=\det((T_1\rightarrow T_2) \rightarrow R)$. By the first assertion, $
\det((T_1\rightarrow T_2) \rightarrow R) = \det(T_1\rightarrow T_2)\cdot 
\det(R)$. Then, $\det(T_1\rightarrow T_2) = \det(T_1^{+}\rightarrow T_2^{
+})$, because $\det(R)=1$. Furthermore, $\det(T_1^{+}\rightarrow T_2^{
+}) = \det(T_1^{+})\cdot\det(T_2^{+})$ by the first assertion.\qed
\end{proof}
   	
	As we have seen above, the converse and switching operations preserve 
unimodularity. Together with the join, these operations generate a 
subclass $\mathcal{H}$ of unimodular tournaments, defined as follows.

\begin{enumerate}
	\item The unique 2-tournament is in $\mathcal{H}$.
  \item If $T_1,T_2$ are in $\mathcal{H}$, then the tournament $T_1
	\rightarrow T_2$ is in $\mathcal{H}$.
  \item The switch of a tournament in $\mathcal{H}$ is also in $\mathcal{H}$.
\end{enumerate}

  Let $T$ be a tournament with $n \geq 4$ vertices. We say that $T$ is 
\emph{switching decomposable} if there exist two tournaments $T_1$ and $T
_2$, each with at least $2$ vertices, such that $T$ is switching 
equivalent to $T_1\rightarrow T_2$. Otherwise, we say that $T$ is
\emph{switching indecomposable}. Switching decomposability coincides with the 
\emph{bijoin decomposability} \cite{bankoussou2017matrix, bui2007unifying}.
  
\begin{example}
For an odd integer $n$, consider the well-known \emph{circular tournament}
$C_n$ whose vertex set is the additive group $\mathbb{Z}_n=\{0,1,\cdots,n
-1\}$ of integers modulo $n$, such that $i$ dominates $j$ if and only if 
$i-j\in\{1,\cdots,(n-1)/2\}$. The tournament $C_n$ is strongly connected.
However, by reversing the arcs between the
even and the odd vertices we obtain a transitive tournament. Hence, $C_n$
is switching decomposable for every odd integer $n\geq5$.
\end{example}

Clearly, every tournament in $\mathcal{H}$ with more than $2$ 
vertices is switching decomposable. It is easy to check that all 
unimodular tournaments with at most $6$ vertices are in $\mathcal{H}$. 
However, we have found a switching indecomposable unimodular tournament 
with $8$ vertices. Its skew-adjacency matrix is $\mathbb{F}_q$
\[
\begin{pmatrix*}[r]
0 & -1 & -1 & -1 & 1 & -1 & -1 & -1 \\
1 & 0 & -1 & -1 & -1 & -1 & -1 & -1 \\
1 & 1 & 0 & -1 & -1 & 1 & -1 & -1 \\
1 & 1 & 1 & 0 & -1 & -1 & -1 & -1 \\
-1 & 1 & 1 & 1 & 0 & -1 & -1 & -1 \\
1 & 1 & -1 & 1 & 1 & 0 & 1 & -1 \\
1 & 1 & 1 & 1 & 1 & -1 & 0 & -1 \\
1 & 1 & 1 & 1 & 1 & 1 & 1 & 0
\end{pmatrix*}
\]

\begin{problem}\label{question:1}
  Find an infinite family of unimodular switching indecomposable
tournaments.
\end{problem}

\section{Spectral properties of unimodular tournaments}
	Let $S$ be an integral skew-symmetric matrix. The nonzero eigenvalues 
of $S$ are purely imaginary and occur as conjugate pairs $\pm i\lambda_1
, \ldots, \pm i\lambda_k$, where $\lambda_1, \ldots, \lambda_k$ are 
totally real algebraic integers. Moreover, the norm $N(\lambda_i)$ of $
\lambda_i$ divides the determinant of $S$.

If $S$ is unimodular, then $\Pi_{i=1}^{k}\lambda_i=\pm1$ and hence, $
\lambda_i$ are all algebraic units. Conversely, suppose  that every $
\lambda_i$ is an algebraic unit. Then, the constant coefficient in the 
minimal polynomial of every eigenvalue is $\pm1$. Hence, the determinant 
of $S$ is $\pm1$. In particular, we have the following result.

\begin{proposition}\label{main10}
	A tournament is unimodular if and only if all its eigenvalues are algebraic units.
\end{proposition}
  	
Godsil \cite{godsil1982eigenvalues} proved that every algebraic integer $
\lambda$ occurs as an eigenvalue of the adjacency  matrix of a digraph. 
Estes \cite{estes1992eigenvalues} proved that if $\lambda$ is a totally 
real integer, that is, all its conjugates are real, then it is an 
eigenvalue of the adjacency matrix of a graph. Recently, Salez
\cite{salez2015every} proved that the graph may be chosen to be a tree.
For tournaments, we can ask the following question.

 \begin{question}
	Let $\lambda$ be a totally real algebraic integer with an odd norm. 
Are there any other conditions on $\lambda$ so that $i\lambda$ is the 
eigenvalue of a tournament?

\end{question}

	A similar question can also be asked about the determinant of 
tournaments. By Proposition \ref{mccarthy}, the determinant of a 
tournament with an even number of vertices is the square of an odd 
number. Conversely,

\begin{question}\label{quest:2}
	Does there exist a tournament whose determinant is $m^{2}$ for every odd
number $m$?
\end{question}
    
	By Theorem \ref{main7}, it is enough to consider Question \ref{quest:2}
for odd prime numbers.

\section{Invertible tournaments}
  Let $T$ be a tournament with skew-adjacency matrix $S$. Let $\phi_{S}(x
)=x^{n}+\sigma_{1}x^{n-1}+\sigma_{2}x^{n-2}+\cdots+\sigma_{n-1}x+\sigma_{
n}$ be the characteristic polynomial of $S$. Then,
\begin{equation}\label{eq:minors}
    \sigma_{k}=(-1)^{k}\sum(\text{all }k\times k\text{ principal minors})\mbox{.}
\end{equation}
  Since $S$ is the skew-adjacency matrix of a tournament, we have
	
\begin{enumerate}
  \item $\sigma_2 = \binom{n}{2}$.
  \item $\sigma_n = \det(S)$.
  \item $\sigma_k = 0$ if $k$ is odd.
\end{enumerate}

	The determinant of a $4$-tournament is $9$ if it is a diamond and $1$ 
otherwise. It follows that $\sigma_4 = 8\delta_T + \binom{n}{4}$, where $
\delta_T$ is the number of diamonds in $T$. In particular, $T$ has no 
diamonds if and only if $\sigma_4 = \binom{n}{4}$.

If $T$ is unimodular, then the inverse $S^{-1}$ of $S$ is an integral 
unimodular skew-symmetric matrix. Furthermore, $\phi_{S^{-1}}(x) = x^{n} + 
\sigma_{n-2}x^{n-2} + \cdots+\sigma_{2}x^{2} + 1$. Hence, a necessary 
condition for $T$ to be invertible is $\sigma_{n-2} = \binom{n}{2}$. The 
following proposition shows that this condition is sufficient.

\begin{proposition}
	Let $T$ be a unimodular $n$-tournament, and let $S$ be its skew-
adjacency matrix. Then, the off-diagonal entries of $S^{-1}$ are odd. 
Moreover, the following assertions are equivalent. 
	\begin{enumerate}[i)]
		\item $T$ is invertible.
		\item Every $(n-2)$-subtournament of $T$ is unimodular.
		\item The coefficient of $x^{2}$ in the characteristic polynomial of $S$ is $\binom{n}{2}$.
	\end{enumerate}
\end{proposition}
  
\begin{proof}
	Let $[n]=\{1,\ldots,n\}$, and let $I$ be a subset of $[n]$. Denote by $
S[I]$ the submatrix of $S$ whose rows and columns are indexed by $I$. 
Let $i\neq j\in [n]$, it follows from Jacobi's complementary minors 
theorem that
	\begin{equation}\label{eq:jacobi}
    \det(S^{-1}[\{i, j\}]) = \det(S[[n]\setminus{\{i, j\}}])\mbox{.}
  \end{equation}
	Moreover, as $S^{-1}$ is skew-symmetric, $\det(S^{-1}[\{i, j\}]) =
(S_{i,j}^{-1})^2$. Then
	\begin{equation}\label{eq:jacobi2}
    (S_{i,j}^{-1})^2 = \det(S[[n]\setminus{\{i, j\}}])\mbox{.}
  \end{equation}
	By Proposition \ref{mccarthy}, $\det(S[[n]\setminus{\{i, j\}}])$ is the square 
of an odd number. Hence, the $(i, j)$-entry of $S^{-1}$ is odd.

 The equivalence $i)\Leftrightarrow ii)$ follows directly from
\eqref{eq:jacobi2}. The equivalence $ii)\Leftrightarrow iii)$ follows from
\eqref{eq:minors} and Proposition \ref{mccarthy} which implies that the 
determinant of tournaments with an even number of vertices is at least $1$.
\end{proof}
  
\begin{example}
  Let $T$ be an $n$-tournament without diamonds and let $S$ be its
skew-adjacency matrix. Every subtournament of $T$ with an even number of 
vertices is unimodular. Hence, $T$ is invertible and $\phi_S(x) = \phi_{S
^{-1}}(x)=x^{n} + \binom{n}{2}x^{n-2} + \cdots + \binom{n}{n-2}x^{2} + 1$.
It follows that $S^{-1}$ is the skew-adjacency matrix of a tournament 
without diamonds.
\end{example}
  
  In the example above, the characteristic polynomial is palindromic, 
that is, the coefficients of $x^{i}$ and $x^{n-i}$ are equal. We call a 
tournament \emph{palindromic} if the characteristic polynomial of its 
skew-adjacency matrix is palindromic. Let $T$ be a tournament and let $S$ be
its skew-adjacency matrix. Clearly, if $\phi_S(x)$ is palindromic,
then $T$ is unimodular, the inverse of $S$ is
the skew-adjacency matrix of a tournament and $\phi_S(x)=\phi_{S^{-1}}({x})$.
In what follows, we give a construction of tournaments whose characteristic 
polynomial is palindromic.
  
  Let $T$ be an $n$-tournament with vertex set $V = \{v_1,\ldots, v_n\}$,
and let $S$ be its skew-adjacency matrix. Let $\hat{T}$ be the 
tournament obtained from $T$ by  adding a copy $T^{\prime}$ of $T$ with 
vertex set $\{v^{\prime}_1,\ldots,v^{\prime}_n\}$, such that $v_{i}$ 
dominates $v^{\prime}_{i}$, and $v_{i}$ dominates $v^{\prime}_{j}$ if 
and only $v_{i}$ dominates $v_{j}$. The skew-adjacency matrix $\hat{S}$ 
of $\hat{T}$ can be written as follows.

\[\hat{S} = 
 \begin{pmatrix}
  S & S + I_n \\
  S - I_n & S
 \end{pmatrix}\mbox{.}
\]

The inverse of $\hat{S}$ is
$
 \begin{pmatrix}
  S & -(S + I_n) \\
  -(S - I_n) & S
 \end{pmatrix}
$
then $\hat{T}$ is invertible. Moreover, $\hat{T}$ and $\hat{T}^{-1}$ are 
switching equivalent. Indeed,
\begin{equation}\label{eq:switch}
\hat{S}^{-1} = D\hat{S}D\mbox{.}
\end{equation}
where
$D=\begin{pmatrix}
    I_n & 0 \\
    0   & -I_n
   \end{pmatrix}
$. It follows that $\phi_{\hat{S}}(x) = \phi_{(\hat{S})^{-1}}(x)$, and
hence $\phi_{\hat{S}}(x)$ is palindromic.

\begin{remark}\label{eq:eq-minors}
  Let $I$ be a nonempty proper subset of $[2n]$. By \eqref{eq:switch},
$\det(\hat{S}[I]) = \det(\hat{S}^{-1}[I])$. Moreover, using Jacobi's 
complementary minors theorem, $\det(\hat{S}^{-1}[I]) = \det(\hat{S}[[2n]
\setminus I])$. It follows that $\det(\hat{S}[I]) = \det(\hat{S}[[2n]
\setminus I])$.
\end{remark}

\section{Embedding of tournaments in unimodular tournaments}
  In the previous section, we proved that every $n$-tournament can be 
embedded in a unimodular $2n$-tournament. For a tournament $T$ on $n$ 
vertices, let $u^{+}(T)$ be the smallest number of vertices we must add 
to $T$ to obtain a unimodular tournament. A dual notion of $u^{+}(T)$ is 
to consider the minimum number $u^{-}(T)$ of vertices we must remove 
from $T$ to obtain a unimodular tournament. It follows from Theorem
\ref{main6} that if $T_1$ and $T_2$ are two tournaments, then
\begin{align}
  u^{+}(T_1\rightarrow T_2) \leq u^{+}(T_1) + u^{+}(T_2)\mbox{,} \\
  u^{-}(T_1\rightarrow T_2) \leq u^{-}(T_1) + u^{-}(T_2)\mbox{.}
\end{align}

  It is shown in \cite{erdos1964representation} that every
$n$-tournament $T$ contains a transitive subtournament of order at least $
\lfloor\log_2(n)\rfloor+1$. In particular, it contains a unimodular 
tournament of order at least $\lfloor\log_2(n)\rfloor$. Then, $u^{-}(T) \leq 
n - \lfloor\log_2(n)\rfloor$.

  The following proposition provides a relationship between $u^{+}(T)$ and $u^{-}(T)$.
\begin{theorem}\label{eq:relation}
  Let $T$ be an $n$-tournament. Then, \[ u^{+}(T) \leq u^{-}(T)\mbox{.} \]
In particular, $u^{+}(T) \leq  n - \lfloor\log_2(n)\rfloor$.
\end{theorem}
  
\begin{proof}
  Let $V = \{v_1,\ldots,v_n\}$ be the vertex set of $T$. Consider the $2n
$-tournament $\hat{T}$ obtained from $T$ and a copy $T^{\prime}$ of $T$ 
as described in the previous section. There exists $I\subset V^{\prime}$
, $|I|=u^{-}(T^{\prime})$, such that $T[V^{\prime}\setminus I]$ is  
unimodular. By Remark \ref{eq:eq-minors}, $\det(\hat{T}[V^{\prime}
\setminus I]) = \det(\hat{T}[V\cup I]) = 1$. Moreover, the tournament $
\hat{T}[V\cup I]$ contains $T$, hence $u^{+}(T) \leq  u^{-}(T)$.
\end{proof}
  
\begin{remark}
  Equality in Theorem \ref{eq:relation} may be strict. Indeed, let $T$ be
the tournament whose skew-adjacency matrix is $S$.
	\[S = 
\begin{pmatrix*}[r]
0 & -1 & -1 & -1 & -1 & -1 & -1 & 1 & -1 \\
1 & 0 & -1 & -1 & -1 & -1 & -1 & -1 & 1 \\
1 & 1 & 0 & -1 & -1 & -1 & 1 & -1 & -1 \\
1 & 1 & 1 & 0 & -1 & -1 & -1 & 1 & -1 \\
1 & 1 & 1 & 1 & 0 & -1 & -1 & -1 & -1 \\
1 & 1 & 1 & 1 & 1 & 0 & -1 & -1 & 1 \\
1 & 1 & -1 & 1 & 1 & 1 & 0 & -1 & -1 \\
-1 & 1 & 1 & -1 & 1 & 1 & 1 & 0 & -1 \\
1 & -1 & 1 & 1 & 1 & -1 & 1 & 1 & 0
\end{pmatrix*}\mbox{.}
\]
By adding a vertex dominating $T$ we obtain a unimodular tournament, hence
$u^{+}(T) = 1$. The tournament $T$ has no unimodular $(n-1)$-subtournament.
Moreover, removing the last three rows of $S$ and their 
corresponding columns yields the skew-adjacency matrix of a unimodular 
tournament, hence $u^{-}(T) = 3$. This example was found using SageMath
\cite{SageMath}.
\end{remark}
  
In what follows, we give a lower bound on $u^{+}(T)$, using the spectra 
of the skew-adjacency matrix of $T$.

\begin{theorem}\label{main8}
  Let $T$ be a non-unimodular $n$-tournament and let $\nu(T)$ be the 
maximum multiplicity among the non-unit eigenvalues of its skew-adjacency
matrix. Then, \[ \nu(T) \leq u^{+}(T)\mbox{.} \]
\end{theorem}
 
  To prove this theorem, we need the following lemma, which is a direct 
consequence of Cauchy Interlace Theorem.

\begin{lemma}\label{main9}
  Let $A$  be a hermitian matrix of order $m$, and let $B$ be a 
principal submatrix of $A$ of order $n$, with an eigenvalue $\lambda$ of 
multiplicity $r$. If $m-n < r$, then $\lambda$ is an eigenvalue of $A$.
\end{lemma}

\begin{proof}[Proof of Theorem \ref{main8}]
  Let $T$ be a non-unimodular $n$-tournament and let $i\lambda$ be a
non-unit eigenvalue of its skew-adjacency matrix $S$ with multiplicity
$\nu(T)$. Let $T^{\prime}$ be an $m$-tournament containing $T$ such that
$m < n + \nu(T)$ and denote by $S^{\prime}$ its skew-adjacency matrix. 
Clearly, $\lambda$ is an eigenvalue of $iS$ with multiplicity $\nu(T)$. 
Then, by Lemma \ref{main9}, $\lambda$ is also an eigenvalue of $iS^{
\prime}$. Hence, $S^{\prime}$ has a non-unit eigenvalue. It follows from 
Proposition $\ref{main10}$ that $T^{\prime}$ is not unimodular. \qed
\end{proof}
  
  Tournaments with large $\nu(T)$ can be obtained from skew-conference 
matrices. Let $T$ be an $n$-tournament and let $S$ be its skew-adjacency 
matrix. Assume that $S$ is a skew-conference matrix.
It follows that the eigenvalues of $S$ are $\pm i\sqrt{n-1}$ each 
with multiplicity $n/2$. As $i\sqrt{n-1}$ is not an algebraic unit, then 
$\nu(T) = n/2$. Hence, by Theorem \ref{main8}, $u^{+}(T)\geq n/2$.
  
  It is conjectured that skew-conference matrices exist if and only if $n
=2$ or $n$ is divisible by $4$ \cite{wallis1971some}. If this conjecture 
is true, Lemma \ref{main9} implies that for every integer $n\geq4$, 
there exists an $n$-tournament $T$ such that $u^{+}(T) \geq \frac{n-3}{2
}$. Denote by $u^{+}(n)$ the maximum $u^{+}(T)$ among $n$-tournaments. 
By the forgoing, we have the following theorem. 

\begin{theorem}
  Assuming the existence of skew-conference matrices of every order 
divisible by $4$, we have
  \[ \frac{n-3}{2} \leq u^{+}(n) \leq n - \lfloor\log_2(n)\rfloor\mbox{.} \]
\end{theorem}
  
  Examples of tournaments with a skew-conference matrix can be obtained 
from Paley tournaments. For a prime power $q\equiv3\mod4$, the Paley 
tournament with $q$ vertices is the tournament whose vertex set is the 
Galois field $GF(q)$, such that $x$ dominates $y$ if and only if $x-y$ 
is a nonzero quadratic residue in $GF(q)$. There are many other infinite 
families of skew-conference matrices, see for example
\cite{koukouvinos2008skew}.

\section{Concluding remarks}
The main concern of this paper is the determinant of the skew-adjacency
matrix of tournaments.
A multiplicative formula for the determinant of the join of two 
tournaments was given. This formula provides a new construction of unimodular 
tournaments. Another construction is the blow-up operation, in which 
every vertex of a tournament is replaced by a tournament with two 
vertices. This construction shows that every $n$-tournament can be 
embedded in a $2n$-unimodular tournament for which the inverse of the 
skew-adjacency matrix is also the skew-adjacency matrix of a tournament.
The minimum number of vertices that must be added to a tournament to be unimodular is
considered. We showed that it does not exceed the minimum number of vertices
to be removed to obtain a unimodular tournament, and that it is related
to the multiplicity of its non-unit eigenvalues.

In addition to the problems presented, many other questions and directions
can be considered.
\begin{itemize}
	\item The construction of the class $\mathcal{H}$, considered in Section
\ref{join}, is simple. Nevertheless, this family seems to be rich as it can
be proven, by induction, that the blow-up of every tournament
is in $\mathcal{H}$. Is a positive proportion of the set of unimodular
tournaments in $\mathcal{H}$?

	\item Find examples of invertible tournaments that are not palindromic.

	\item The problem of finding $u^{+}(T)$ seems extremely hard. We suspect
that there is no polynomial time algorithm to solve this problem. Find a
non-brute force algorithm to compute $u^{+}(T)$.

	\item As we have seen above, skew-conference matrices have non-unit 
eigenvalues with maximum possible multiplicities. Another property of 
skew-conference matrices is that they have maximum determinant among
zero-diagonal $\{-1, 1\}$-matrices. Do tournaments with skew-conference
adjacency matrices have maximum $u^{+}(T)$?
\end{itemize}

\bibliographystyle{plain}
\bibliography{biblio}

\end{document}